\documentclass[11 pt, psamsfonts]{amsart}
\usepackage{amssymb}

\def\ignore#1{\relax}

\def\g{\mathfrak g}
\def\h{\mathfrak h}

\def\R{{\mathbb R}}
\def\Z{{\mathbb Z}}

\def\C{{\mathbb C}}

\def\la{\lambda}
\def\La{\Lambda}

\def\N{\mathbb N}

\def\Ca{\mathcal C}
\def\Cl{Cl}

\def\one{\mathbf 1}

\def\b{{\bf b}}

\def\U{{\bf U}}

\def\v{\vskip 2.5mm}

\def\ignore#1{\relax}
\def\one{{\bf 1}}

\def\om{\omega}
\def\e{\epsilon}
\def\1{{\bf 1}}
\def\lan{\langle}
\def\ra{\rangle}

\def\End{{\rm End}}

\def\ve{\varepsilon}

\calclayout

\renewenvironment{proof}[1][\proofname]{\par
  \normalfont
  \topsep6\p@\@plus6\p@ \trivlist
  \item[\hskip\labelsep\bfseries
    #1\@addpunct{.}]\ignorespaces
}{%
  \qed\endtrivlist
}

\def\th@plain{%
  \let\thmhead\thmhead@plain \let\swappedhead\swappedhead@plain
  \thm@preskip.5\baselineskip\@plus.2\baselineskip
                                    \@minus.2\baselineskip
  \thm@postskip\thm@preskip
  \itshape
\renewcommand{\labelenumi}{{(\alph{enumi})\quad}}
                        \renewcommand{\labelenumii}{{(\roman{enumii})\ }}
}
\def\th@definition{%
  \let\thmhead\thmhead@plain \let\swappedhead\swappedhead@plain
  \thm@preskip.5\baselineskip\@plus.2\baselineskip
                                    \@minus.2\baselineskip
  \thm@postskip\thm@preskip
  \upshape
}
\def\th@remark{%
  \thm@headfont{\itshape}
  \let\thmhead\thmhead@plain \let\swappedhead\swappedhead@plain
  \thm@preskip.5\baselineskip\@plus.2\baselineskip
                                    \@minus.2\baselineskip
  \thm@postskip\thm@preskip
  \upshape
}

{\theoremstyle{plain}
\newtheorem{theorem}{Theorem}[section]
}

{\theoremstyle{plain}
\newtheorem{proposition}[theorem]{Proposition}
}

{\theoremstyle{plain}
\newtheorem{corollary}[theorem]{Corollary}
}

{\theoremstyle{plain}
\newtheorem{lemma}[theorem]{Lemma}
}

{\theoremstyle{plain}

}

{\theoremstyle{definition}
\newtheorem{definition}[theorem]{Definition}
}

{\theoremstyle{definition}

}

{\theoremstyle{remark}
\newtheorem{remark}[theorem]{Remark}
}

{\theoremstyle{remark}

}

\numberwithin{equation}{section}

\renewcommand{\labelenumi}{{ \theenumi.}}
\renewcommand{\labelenumii}{{(\alph{enumii})}}

\def\la{\lambda}
\def\al{\alpha}

\def\so{{\mathfrak s}{\mathfrak o}}

\def\choose #1 #2{\begin{pmatrix}#1\\#2\end{pmatrix}}

\def\lad{{\la^{\dagger}}}
\def\mud{{\mu^{\dagger}}}
\def\y{{\bf y}}
\def\x{{\bf x}}

\ignore{
\input BoxedEPS
\SetRokickiEPSFSpecial
\HideDisplacementBoxes
\SetEPSFDirectory{./graphics/} }

\ignore{
\input BoxedEPS
\SetTexturesEPSFSpecial
\HideDisplacementBoxes
\SetEPSFDirectory{:graphics:}
}

\def\N{\mathbb N}

\def\Ca{\mathcal C}
\def\Cl{Cl}

\def\one{\mathbf 1}

\def\b{{\bf b}}

\def\x{{\bf x}}
\def\y{{\bf y}}

\def\U{{\bf U}}

\def\v{{\bf v}}

\def\ignore#1{\relax}
\def\one{{\bf 1}}

\def\om{\omega}
\def\e{\epsilon}
\def\1{{\bf 1}}
\def\lan{\langle}
\def\ra{\rangle}

\def\End{{\rm End}}

\def\ve{\varepsilon}

\def\SS{\tilde S}

\begin{document}

\title[Centralizer Algebras for Spin Representations ]
{On Centralizer Algebras for Spin Representations }

\author{Hans Wenzl}
\thanks{$^*$Supported in part by NSF grants.}

\address{Department of Mathematics\\ University of California\\ San Diego,
California}

\email{hwenzl@ucsd.edu}

\begin{abstract}
We give a presentation of the centralizer algebras for 
tensor products of spinor
representations of quantum groups
via generators and relations. In the even-dimensional 
case, this can be described in terms of non-standard
$q$-deformations of orthogonal
Lie algebras; in the odd-dimensional case only a certain subalgebra
will appear. In the classical case $q=1$ the relations boil
down to Lie algebra relations.
\end{abstract}
\maketitle

Classically, representations of Lie groups were studied by decomposing
tensor products of a simple generating representation. This worked 
very well for  the vector representations of the general  linear group
and  to a somewhat lesser extent also for  the vector representations
of orthogonal and symplectic groups. More recently, $q$-versions of 
these centralizer algebras
were studied in connection with Drinfeld-Jimbo quantum groups
which have had applications in a number of fields.
In this paper, we study centralizer algebras for spinor representations.
Some of the possible applications will be given below.

We first give a description of the centralizer algebras for
the classical spinor groups. If $N$ is even,
this can be comparatively easily deduced from results
by Hasegawa \cite{Has} for the action of $Pin(N)$.
It follows that the commutant of its action on
the $l$-th tensor power $S^{\otimes l}$ of its spinor
representation $S$ is given by a representation of $so_l$.
Moreover, we give a precise identification of the action
of generators which is compatible with embeddings
$\End(S^{\otimes l})\subset \End(S^{\otimes ( l+1)})$.
This is important for studying the corresponding tensor
categories and is not immediately obvious from
the results in \cite{Has}.
We also prove an analogous result for 
the odd-dimensional case, which is more complicated.
Here the
commutant is generated by a subalgebra of the universal
enveloping algebra $Uso_l$.

We then extend these results to the setting of quantum groups.
In principle, it should be possible to do this similar to the classical
case, using known $q$-deformations
of Clifford algebras, see Section \ref{qCliffapp}.
However, due to their complicated multiplicative
structure we determine generating
elements  via a straightforward approach.
They are $q$-deformations of the canonical element
$\sum e_i\otimes e_i\subset Cl(V)\otimes Cl(V)\cong \End(S^{\otimes 2})$,
where $(e_i)_i$ is an orthonormal basis of the vector representation
$V$ and $Cl(V)$ is its Clifford algebra.
We obtain from this elements which satisfy the relations of
generators of another $q$-deformation
$U_q'so_l$ of the  universal enveloping algebra $Uso_l$. It has appeared before
in work of Gavrilik and Klimyk \cite{GK}, Noumi and Sugitani \cite{NS} and Letzter \cite{Le}.
Unlike the usual $q$-deformation of $Uso_l$, it does not have
a Hopf algebra structure. Our main result is that we again have
a duality between the actions of $U_qso_N\rtimes \Z/2$
and $U'_qso_l$ acting as each others commutants on
$S^{\otimes l}$ for $N$ is even. 
Again, the situation is  more complicated in the odd-dimensional case
where we have to consider a subalgebra of $U'_qso_l$.

It is worthwhile mentioning that one of the problems
for spinor representations is that already their second
tensor power contains an increasing number of 
irreducible representations. This makes it difficult
to characterize the centralizer algebras via
braid representations, which worked well for vector
representations. A similar problem was encountered
by Rowell and Wang in their study of certain braid
representations, which they conjecture to be
related to spinor representations at certain roots
of unity, see \cite{RW}. Our results should be useful
in studying this question. This and other potential
applications are discussed at the end of this paper.

Here is the paper in more detail: We first  show how the
commutant on the $l$-th tensor product of a spinor representation
is related to $so_l$ by fairly elementary methods. While
many of the results are not new, it serves as a blueprint
for the more difficult quantum group case.
We then review  basic material 
from the study of Lie algebras and Drinfeld-Jimbo quantum groups.
This is then used to prove the already sketched duality results 
for quantum groups, where we find generators for
$\End_\U(S^{\otimes l})$ in the third section, and relations
in the fourth section.

$Acknowledgements:$ Part of the work on this paper was done
while the author was visiting Aarhus University and the
Hausdorff Institute. The author is grateful to these institutions
for hospitality and support. He would also like to thank
A.J. Wassermann, 
H. H. Andersen and G. Lehrer for references. Thanks are also
going to the referee and to Eric Rowell for pointing out confusing
inconsistencies in notations in an earlier version.

\section{Duality for Spinor Representations}\label{classicalcase}
We assume throughout this paper all the algebras to be defined over the field
of complex numbers, with $q$ not being a root of unity. For
possible generalizations, see Remark \ref{generalization}.

\subsection{Clifford Algebras and spinor representations}
Let $\{ e_1,\ ...\ e_N\}$ be an orthonormal basis of the finite
dimensional inner product space $V$. Then the Clifford algebra
$Cl =Cl(N)$ corresponding to $V$ can be defined via generators,
also denoted by $e_i$, and relations
$$e_ie_j+e_je_i=2\delta_{ij}1,\quad {\rm for}\ 1\leq i,j\leq N.$$
It is well-known that $Cl(N)$ has dimension $2^N$. It is isomorphic
to $M_{2^{N/2}}$ for $N$ even, and to $M_{2^{(N-1)/2}}\oplus M_{2^{(N-1)/2}}$
for $N$ odd; here $M_d$ denotes the $d\times d$ matrices.
The action of an element $g$ in the orthogonal group $O(N)$
on $V$ induces an automorphism $\al_g$ on $Cl(N)$, for each $g\in O(N)$.
As any automorphism of the $d\times d$ matrices is inner,
we obtain a projective representation of $O(N)$ on a
$2^{N/2}$ dimensional module $S$ in the {\it even-dimensional case}.
By restriction, the module $S$ becomes a projective $O(N-1)$-module
which we will denote by $\tilde S$. It decomposes into the direct sum
of two simple projective $O(N-1)$-modules $\tilde S_+\oplus \tilde S_-$.
These two modules are isomorphic as projective $SO(N-1)$-modules
and simple.

A more direct way to construct the simple projective $SO(N)$ module
for $N$ odd goes as follows: We replace
the full Clifford algebra, which is not simple in the odd-dimensional case, by
$\Cl_{ev}(N)$, the span of all products of $e_i$'s with an even
number of factors. Observe that for $N$ odd, the element $f_N
=e_1e_2\ ...\ e_N$
is in the center of $\Cl(N)$, and $(\gamma f_N)^2=1$ for a suitable
$\gamma\in \{ 1,\sqrt{-1}\}$.  Hence the map induced by
$e_i\mapsto \gamma e_if_N$ defines a homomorphism from $\Cl(N)$ into
$\Cl_{ev}(N)$. If one restricts this homomorphism to the 
simple subalgebra $M_{2^{(N-1)/2}}\cong\Cl(N-1)\subset \Cl(N)$, 
it is obviously not the zero map. Hence
it becomes an isomorphism, by dimension count and simplicity of
$\Cl(N-1)$.
The automorphism $\al_g$ induced by an element $g\in O(N)$
restricts to the subalgebra $\Cl_{ev}(N)$ of  $\Cl(N)$.
Hence, again, we obtain  a projective representation of $O(N)$ on an
$2^{(N-1)/2}$ dimensional module $S$; it is isomorphic as a projective
$SO(N)$ module to the modules $\tilde S_{\pm}$ of the last paragraph.

The projective representations just mentioned extend to honest
representations of suitable covering groups. It will be convenient
to consider these modules over $Pin(N)$ for $N$ even, and over
$Spin(N)$ for $N$ odd. Here $Pin(N)$ and $Spin(N)$ are the
two-fold covering groups of $O(N)$ and $SO(N)$ respectively.

\subsection{Hasegawa's results} Let now $V=\C^{Nl}=\C^N\otimes \C^l$.
Then we have obviously commuting actions of $O(N)$ and $O(l)$ on
$V$, acting on the corresponding tensor factors.

\begin{theorem}\label{Hasegawa} (Hasegawa \cite{Has}) The algebras generated
by the actions of $O(N)$ and $O(l)$ on $Cl(Nl)$ are each others
commutant.
\end{theorem}

It is easy to check that $Cl(Nl)\cong Cl(N)^{\otimes l}
\cong Cl(l)^{\otimes N}$ as vector spaces. This strongly suggests
a relationship between the commutant of the action of $Spin(N)$
on the $l$-fold tensor product of its spinor representation,
and the group $Spin(l)$.  Observations to this extent have been
made at the combinatorial level in several papers before, e.g. \cite{Ba}, \cite{BS}.
However, the precise result we need is a little bit more subtle
and does not immediately follow from the results above.
In particular, in our context there are nontrivial distinctions between
the odd and even-dimensional cases which do not occur in \cite{Has}.

\subsection{Some elementary lemmas} Let $f_m=e_1e_2\ ...\ e_m\in Cl(m)$
for $m\in\N$. Moreover, consider the map
\begin{equation}\label{Phidef}
\Phi : 1\otimes \ ...\ \otimes e_i\otimes\ ...\ \otimes 1\in Cl(N)^{\otimes l}
\quad \mapsto \quad
\begin{cases}  f_{(j-1)N}e_{(j-1)N+i} & \text{if $j$ is odd,}\\
f_{jN}e_{(j-1)N+i} &\text{if $j$ is even,}
\end{cases}
\end{equation}
where $e_i$ is in the $j$-th factor of $Cl(N)^{\otimes l}$.
Then we have the following easy lemma:

\begin{lemma} \label{embedding} (a) $f_me_i=(-1)^me_if_m$ for $i>m$, and  
$f_me_i=-(-1)^me_if_m$ for $i\leq m$.

(b) If $N$ is even, the map $\Phi$ above extends to an algebra 
and $O(N)$-module  isomorphism
between $Cl(N)^{\otimes l}$ and $Cl(Nl)$, and
 $\Phi(f_N\otimes f_N)=(-1)^{N(N-1)/2}f_{2N}$.

(c) For $N$ odd, the map $\Phi$ defines an embedding of 
$\Cl_{ev}(N)^{\otimes l}$ into  $\Cl(Nl)$. It maps the element
$e_re_s$, in the $j$-th factor of $\Cl_{ev}(N)^{\otimes l}$,
to $e_{(j-1)N+r}e_{(j-1)N+s}$.
\end{lemma}

$Proof.$ Part (a) is straightforward. For part (b), one first
checks that the map $\Phi$ indeed defines a nonzero algebra homomorphism.
This is straightforward. As both source and target algebras are simple
and have the same dimension, $\Phi$ is an isomorphism.
As $\Phi$ is an $O(N)$ module morphism on the linear
span of the generators, and $O(N)$ acts via algebra automorphisms,
$\Phi$ is also an $O(N)$ morphism.
The second part of (b) is checked easily using (a) 
and the definition of $\Phi$,
\ref{Phidef}. Part (c) again is straightforward.

\medskip

Recall that the Lie algebra $so_l$ is isomorphic to the subset of
$l\times l$ matrices spanned by $L_{rs} = E_{rs}-E_{sr}$, $1\leq r<s\leq l$,
where the $E_{rs}$ are matrix units. It can also be defined via generators
$L_1,L_2,\ ...\ L_{l-1}$ and relations $[L_i,[L_i,L_{i\pm 1}]]=-L_{\pm 1}$
and $[L_i,L_j]=0$ for $|i-j|>1$. Indeed, it is easy to check that these
relations are satisfied for $L_i=L_{i,i+1}$. Also observe that
one can replace $-L_{i\pm 1}$ by $L_{i\pm 1}$ on the right hand side
of the first relation after substituting $L_i$ by $\sqrt{-1}L_i$.
Let now $N=2k$ be $even$.
We  define the elements $C_{rs}\in Cl(Nl)$ by
$C_{rs} = \frac{1}{2}\sum_i^N e_{(r-1)N+i} e_{(s-1)N+i}$ and
$C'_{rs} = \frac{1}{2}\sum_i^{N-1} e_{(r-1)N+i} e_{(s-1)N+i}$
for $1\leq r<s\leq l$.
Then we have

\begin{lemma}\label{solrelations}
The elements $C_{rs}$ and $C'_{r,s}$ satisfy the commutation relations of the generators
of the Lie algebra $so_l$.
\end{lemma}

$Proof.$ If for indices $p,q,r,s$ the
set $\{ r,s\}\cap\{p,q\}$ is empty or has two elements
$[L_{rs},L_{pq}]=0$. Otherwise, if, say $s=p$, we
get $[L_{rs},L_{sq}]=L_{rq}$. But then we also have
$[C_{rs},C_{sq}]= $
$$=\frac{1}{4}\sum_{i\neq j} 
e_{(r-1)N+i}e_{(s-1)N+i}  e_{(s-1)N+j}e_{(q-1)N+j} 
-
e_{(s-1)N+j}e_{(q-1)N+j} e_{(r-1)N+i} e_{(s-1)N+i}$$ 
$$+ \frac{1}{4}\sum_{i} 2 e_{(r-1)N+i}e_{(q-1)N+j}.$$
One checks that the first sum is equal to 0, and the second
one is equal to $C_{rq}$, which is the required relation. The proof for
the $C'_{rs}$ goes the same way.

\medskip
We shall need the
precise preimages of $C_{12}$ and $C_{23}$ under the isomorphism $\Phi$.
It follows from Lemma \ref{embedding} that they are given by
$$\Phi^{-1}(C_{12})=(-1)^{N(N-1)/2}\sum_{i=1}^N e_if_N\otimes e_if_N\otimes 1
\quad {\rm and}
\quad \Phi^{-1}(C_{23})=\sum_{i=1}^N 1\otimes e_i\otimes e_i.$$
Similarly, the elements $\Phi^{-1}(C'_{12})$ and  $\Phi^{-1}(C'_{23})$
are given by the same sums, now only going until $N-1$.

\begin{corollary}\label{Ccorollary}
The elements $\Phi^{-1}(C_{12})$ and $\Phi^{-1}(C_{23})$
 are in $\End_{Pin(N)}(S^{\otimes 3})$, and the elements
$\Phi^{-1}(C'_{12})$ and $\Phi^{-1}(C'_{23})$ are  in 
 $\End_{Pin(N-1)}(\tilde S^{\otimes 3})$, where $\tilde S$
is $S$ viewed as a $Pin(N-1)$-module.
\end{corollary}

$Proof.$ As $g\in O(N)$ fixes $\sum e_i\otimes e_i$, viewed as an element
in $V^{\otimes 2}$, one deduces that conjugation of
$\Phi^{-1}(C_{23})$ by a lift of $g$ in $Pin(N)$ leaves it invariant.
The other statements follow similarly.

\ignore{
$Proof.$ As any two elements $e_ie_j$ and $e_re_s$ in $\Cl(kl)$ commute
whenever $\{ i,j\}\cap \{ r,s\}=\emptyset$, the map $\Phi$
does indeed define an algebra homomorphism from $\Cl_{ev}(k)^{\otimes l}$
into  $\Cl(kl)$ which is injective.
It follows from the discussion before the lemma that
$\End_{Spin(k)}^{\otimes l}\cong \Cl_{ev}(k)^{\otimes l}$.
}

\subsection{Eigenvalues of $C$} It remains to determine
the structure of the representation of $so_l$ in $Cl(Nl)$.
For this, we define the elements $C$, $C'$ and $\tilde C_m$ 
in $Cl(N)^{\otimes 2}$  by
\begin{equation}\label{Cdefi}
C\ =\ \frac{1}{2}\ \sum e_i\otimes\ e_i, 
\end{equation}
\begin{equation}\label{Cmdef}
\tilde C_m=m!\ \sum_{i_1<i_2<\ ...\ <i_m} e_{i_1} e_{i_2}\ ...\  e_{i_m}\otimes
 e_{i_1} e_{i_2}\ ...\  e_{i_m},
\end{equation}
and $C'$ is defined as $C$, but with the summation only going until $N-1$.
Observe that $\tilde C_1=2C$. Moreover, observe that 
the elements $e_i\otimes e_i$ and  $e_j\otimes e_j$
commute also for $i\neq j$.
We will also need the polynomials
$P_m(N,x)$ defined inductively by $P_0(N,x)=1$, $P_1(N,x)=x$
and
\begin{equation}\label{recursionp}
P_{m+1}(N,x)=xP_m(N,x)+m(N+1-m)P_{m-1}(N,x).
\end{equation}
Then we have

\begin{proposition}\label{eigenvalprop} The element $C\in \Cl(N)^{\otimes 2}$
has the eigenvalues $iN/2,i(N/2-1),\ ...,\ i(1-N/2),-iN/2$.
The same statement holds for the element $C'$, with $N$ replaced 
by $N-1$.
\end{proposition}

$Proof.$ Let us first prove the following recursion relation:
\begin{equation}\label{recursiona}
\tilde C_1\tilde C_m=\tilde C_{m+1}+ m(N+1-m)\tilde C_{m-1}.
\end{equation}
Observe that if we define $y_j=e_j\otimes e_j$, then
$y_jy_i=y_iy_j$ and $y_i^2=1$. Moreover, $\tilde C_m=
m!\ \sum_{i_1<i_2<\ ...\ <i_m\leq N} y_{i_1} y_{i_2}\ ...\  y_{i_m}$.
Now if one multiplies  a monomial $ y_{i_1} y_{i_2}\ ...\  y_{i_m}$
by $y_i$, we will obtain a monomial with $m+ 1$ factors if
$i\neq \{ i_j, 1\leq j\leq m\}$, and one with $m-1$ factors otherwise.
By symmetry $\tilde C_1\tilde C_m$ is a linear combination of 
$\tilde C_{m+1}$ and $\tilde C_{m-1}$.
It remains to calculate the coefficients of the leading terms
$y_1y_2\ ...\ y_{m\pm 1}$, which is easy.

Comparing \ref{recursiona} with  \ref{recursionp}, we see that
$\tilde C_m=P_m(N,\tilde C_1)$ for $m\leq N$. As $\tilde C_{N+1}=0$,
it follows that $P_{N+1}(N,\tilde C_1)=0$. As moreover $\tilde C_1^m$
has degree $m$ as a polynomial in the $y_i$s, 
the elements $\tilde C_1, \tilde C_1^2, ...\ \tilde C_1^N$
are linearly independent. It follows
that $P_{N+1}(N,x)$ is the characteristic polynomial of $\tilde C$.
The statement now follows from Lemma \ref{sl2calculations},
which will be proved in the next subsection. The proof for $C'$ goes 
exactly the same way, with $\tilde C'$ defined as in  (\ref{Cmdef})
with the indices only going until $N-1$.

\subsection{Structure coefficients} This subsection serves
to calculate the eigenvalues of the polynomials $P_{N+1}(N,x)$.
Moreover, we do some additional calculations which are
useful for an explicit description of $\End_{(S)pin(N)}(S^{\otimes l})$
also in the quantum case. All of this is obtained in a fairly straightforward
way from the representation theory of $sl_2$, which is well-known
(see e.g. \cite{Hum}). Presumably, most of the results in this section
are known to experts.

Let $H,E,F$ be the usual generators of the Lie algebra $sl_2$,
and let $V_N$ be its $(N+1)$-dimensional simple representation.
It can be defined via a basis $\{ e_o, e_1,\ ...\ e_N\}$
of eigenvectors of $H$ which satisfies
\begin{equation}\label{sl2rep}
E.w_r=(N-r+1)e_{r-1},\quad H.e_r=(N-2r)e_r,\quad F.e_r=(r+1)e_{r+1}.
\end{equation}
The following lemma is well-known and easy to check
(e.g. (a) follows from the fact that $E-F$ is conjugate to $iH$).

\begin{lemma}\label{sl2facts} (a) The element $E-F$ has eigenvalues
$(N-2r)i$, $0\leq r\leq N$ in the representation $V_N$.

(b) The elements $iH/2$ and $(E\pm F)/2$ satisfy the relations
of the generators $L_{rs}$ of the Lie algebra $so_3$ for $1\leq r<s\leq 3$.
\end{lemma}

\begin{lemma}\label{sl2calculations}
(a) The polynomial $P_{N+1}(N,x)$ has the roots $(N-2r)i$, $0\leq r\leq N$.

(b) Let $\x=(x(\la)_r)$ and $\y=(y(\la)_r)$ be the right 
and left eigenvectors of
$E-F$ for the eigenvalue $\la$, with respect to the basis $(e_r)$
and normalized by $\x_0=1=\y_0$. Then
$$x(\la)_r=\frac{(N-r)!P_r(\la)}{N!}\quad {\rm and}\quad
y(\la)_r=\frac{P_r(\la)}{r!}.$$
\end{lemma}

$Proof.$ Writing $E-F$ as a matrix given by \ref{sl2rep}, we obtain
from $(E-F)\x(\la)=\la\x(\la)$ the recursion relation $x_o=1$, $x_1=\la/N$
and 
$$x_{r+1}=\frac{1}{N-r}(\la x_r+rx_{r-1}).$$
Similarly, one obtains from $\y^t(E-F)=\la\y^t$ the recursion relation
$y_o=1$, $y_1=\la$ and $y_{r+1}=(-\la y+(k+1-r)y_{r-1})/(r+1)$.
Comparing this with the recursion relation \ref{recursionp},
one can easily check claim (b). Moreover, we obtain from the
last coordinate in the equation  $(E-F)\x(\la)=\la\x(\la)$
that $\la x_N(\la)+Nx_{N-1}(\la)=0$. Hence
$$0=(\la P_N(N,\la) + N\cdot 1 P_{N-1}(N,\la))/N! = P_{N+1}(N,\la)/N!$$
for any eigenvalue $\la$ of $E-F$.
This together with Lemma \ref{sl2facts}
implies statement (a). 
\smallskip

\begin{proposition}\label{solrep} (a)
Let $N$ be even. Then we obtain a representation of $so_l$
in $\End_{Pin(N)}S^{\otimes l}$ by mapping the element
$L_{rs}$ to the inverse image of $C_{rs}$ under the
isomorphism $\Phi$. For $l=3$, it contains an irreducible 
$(N+1)$-dimensional representation of $so_3$.

(b) If $N$ is odd, we obtain a representation of the subalgebra
of the universal enveloping algebra $Uso_l$ generated by
the elements $L_{rs}^2$, $1\leq r<s\leq l$, by mapping these
generators to the inverse images of the elements $(C_{rs}')^2$
in $\End_{Spin(N)}S^{\otimes l}$. For $l=3$, it contains an irreducible 
$(N+1)/2$-dimensional representation of this subalgebra.
\end{proposition}

$Proof.$ If $N$ is even, we have $\End(S^{\otimes l})\cong Cl(N)^{\otimes l}$,
and  $\End_{Pin(N)}(S^{\otimes l})\cong Cl(N)^{\otimes l}_{\1}$,
the component of $Cl(N)^{\otimes l}\cong Cl(Nl)$
on which the multiplicative action of $O(N)$ is trivial.
By construction, the elements $C_{rs}$ are fixed by the
action of $O(N)$, and they define a representation of $so_l$
by Lemma \ref{solrelations}. By Prop. \ref{eigenvalprop} and Lemma
\ref{sl2rep},
the largest eigenvalue of the image of $H$ is $N$, which shows the existence
of an irreducible $(N+1)$-dimensional representation of $so_3\cong sl_2$
in $Cl(N)^{\otimes 3}$.

For $N$ odd,  we obtain a representation of $so_l$ in 
$\End_{Pin(N)}(\tilde S^{\otimes l})$ by Lemma \ref{solrelations}
and its corollary. Up to a common sign, the elements 
$\Phi^{-1}( (C_{rs}')^2)$
coincide with the elements $\Phi^{-1}(C_{rs}^2)\in\Cl_{ev}(N)^{\otimes l}
\cong \End(S^{\otimes l})$. 
Moreover, it is
easy to see that $(E-F)^2$ and $H^2$ leave the spans of the even and
of the odd basis vectors invariant, and that these are irreducible submodules.

\begin{definition}\label{algdef} (a) We define the algebra $U(l,k)$
via generators $D_1, D_2,\ ...\ D_{l-1}$ and relations
$[D_i,[D_i,D_{i\pm 1}]]=D_{i\pm 1}$, $[D_i,D_j]=0$ if $|i-j|>1$
and by $\prod_{j=-k}^k (D_i-j)=0$.

(b) We define the algebra $Uo(l,k)$ as a subalgebra, generated by
 $D^2_1, D^2_2,\ ...\ D^2_{l-1}$, of an algebra
with generators $D_1, D_2,\ ...,\ D_{l-1}$, where the $D_i$ satisfy the relations
 $[D_i,[D_i,D_{i\pm 1}]]=D_{i\pm 1}$, $[D_i,D_j]=0$ if $|i-j|>1$
and by $\prod_{j=-k+1}^k (D_i-j-1/2)=0$.
\end{definition}

Observe that the algebra $U(l,k)$ is a quotient of the universal enveloping algebra
$Uso_l$ of the Lie algebra $so_{l}$, while $Uo(l,k)$ is a quotient of a subalgebra of
$Uso_l$.

\begin{corollary}\label{genrelations} (a) If $N=2k$ is even, the image of $Uso_l$ in
$\End_{Pin(2k)}(S^{\otimes l})$ is a quotient of $U(l,k)$.

(b) If $N=2k-1$, the image of $Uso_l$ in
$\End_{Spin(2k-1)}(S^{\otimes l})$ is a quotient of $Uo(l,k)$.
\end{corollary}

\section{Lie algebras and quantum groups} 

\subsection{Quantum groups}
We list some basic information
about quantum groups (see e.g. \cite{jantzen}, \cite{Lu}).
Let $\g$ be a symmetrizable
Kac-Moody algebra given by a Coxeter graph $X$ with $k$ vertices,
with generators $e_i$ and $f_i$, $1\leq i\leq k$;
eventually, we will only be interested in orthogonal Lie algebras.
We denote the simple roots by $\alpha_i, i=1, ..., k$. Fix an
invariant bilinear form $\lan\ ,\ \ra$ on $\h^*$, and define $\check{\alpha}_i
=\frac{2}{\lan\alpha_i ,\alpha_i\ra}\alpha_i$. If all the roots have the same
length, we assume $\lan\ ,\ \ra$ to be normalized such that  $\check{\alpha}_i
=\alpha$.  If $\lan\ ,\ \ra$ is nondegenerate, we define
the fundamental weights $\La_j$ by $\lan\check{\alpha}_i,\La_j\ra
=\delta_{ij}$. If $\al\in\h^*$, the reflection $s_\al$ on $\h^*$ is defined by
$s_\al(\la)=\la-\lan \la,\check{\alpha}\ra\al$.
We denote by $\U = U_q\g$
the Drinfeld-Jimbo quantum group corresponding to the semmisimple
Lie algebra $\g$.  It is well-known that for $q$ not a root of unity,
the representation theory of $\U$ is essentially the same as the one
of $\g$, i.e. same labeling set of simple representations,
character formulas etc.  So we will sometimes  state results only for $\g$
when its generalization to $\U$ is obvious.

\subsection{Gradation via Lie subalgebra} \label{Lieprelim}
Let $\g_0$ be a Lie subalgebra of $\g$ corresponding to the graph
obtained from $X$ by removing the vertex labeled by $1$.
If $\la=\sum_{i=1}^N a_i\La_i$ is a weight of $\g$, we denote
by $\hat\la=\sum_{i=2}^N a_i\La_i$ the corresponding weight of
$\g_0$ (after obvious identifications of the fundamental weights
of $\g_0$ with a subset of fundamental weights of $\g$).
Let $V$ be a finite dimensional module of $\g$ such that
 $\langle \La ,\La_1\rangle =c$, a constant for
all highest weights $\La$ of $V$;
here $\La_1$ is the fundamental weight corresponding to $1$.
We denote by $V[0]$ the $\g_0$-module generated by the highest weight
vectors of $V$. 
More generally, we define the level $i$ subspaces  $V^{\otimes n}[i]$
for tensor powers of $V$ and for $i=0,1,2, ...$ by
$$V^{\otimes n}[i]\ =\ span\{ V^{\otimes n}[\mu ],
\ \lan n\La-\mu,\La_1\ra=i\}.$$
It is easy to see that $V^{\otimes n}[0]=(V[0])^{\otimes n}$ for all
$n\in\N$.
Conversely, we say that a weight $\mu$ has level $i$ in $V^{\otimes n}$
if $V^{\otimes n}[\mu]\subset V^{\otimes n}[i]$; in this case we
denote the level of $\mu$ by $lev_n(\mu)$ or just $lev(\mu)$ if
no confusion arises.

\begin{lemma}\label{prelim} Let $V$ be a $\g$-module
as just described. Then

(a) $V^{\otimes n}[0]\cong V_{\hat\La}^{\otimes n}$ as a $\g_0$-module.

(b) Let $W$ be the $\g$-module generated by  $V^{\otimes n}[0]$.
Then  $\End_\g(W)\cong \End_{\g_0}(V_{\hat\La}^{\otimes n}[0])$.
In particular, $mult_{V_\mu}(V^{\otimes n})=
mult_{V_{\hat{\mu}}}(V_{\hat{\La}}^{\otimes n})[0]$ for any
weight $\mu$ with $lev_n(\mu)=0$.
\end{lemma}

$Proof.$ Part (a) was already shown in \cite{exc}, Lemma 1.1 if
$V$ is irreducible. The proof carries over easily to our slightly
more general setting.
Part (b) is an easy consequence of part (a).

\subsection{Traces and contractions}\label{tracecont}  The following material can be
found in e.g. \cite{Ks}, \cite{Turaev} and \cite{OW}, Section 1.4.
Let $W$ be  a $\U$-module, and let $a\in\End_U(W)$. Then 
the categorical trace or $q$-trace $Tr_q(a)$ is given by
$Tr_q(a)=Tr(q^\rho a)$; here $Tr$ is the ordinary trace on
$\End(W)$, and $q^\rho$ acts on the weight vector $w\in W$ with weight
$\mu$ by the scalar $q^{\langle \mu, \rho\rangle}$.
Let  $W=V_\la$ be an irreducible module with highest weight $\la$.
Using the notation $[n]=(q^n-q^{-n})/(q-q^{-1})$, we can explicitly
write the $q$-dimension as
\begin{equation}\label{qdim}
\dim_qV_\la = \prod_{\al >0}\frac{[\langle \la+\rho, \al\rangle]}
{[\langle \rho,\al\rangle]}.
\end{equation}
In particular, if $e$ is a minimal idempotent in $\End_U(W)$
projecting onto an irreducible submodule $\cong V_\la$ of $W$,
we have $Tr_q(e)=\dim_qV_\la$. The normalized trace $tr_q$ is defined by
$tr_q=(1/\dim_q W)Tr_q$.

Let $A\subset B$ be finite-dimensional semisimple algebras with a 
nondegenerate normalized trace $tr$ on $B$ such that also its restriction
to $A$ is nondegenerate; nondegenerate here means that the bilinear form
$\langle b_1,b_2\rangle = tr(b_1b_2)$ is nondegenerate. Then the
orthogonal projection from $B$ onto $A$ with respect to this bilinear
form is usually called the trace preserving conditional expectation
$\ve_A$. Its values are uniquely determined by
$tr(a\ve(b))=tr(ab)$ for all $a\in A$ and all $b\in B$.

In the setting above, one can define an algebra extension $B_1$ of $B$
with respect to the inclusion $A\subset B$, Jones' basic construction,
as follows: It is generated by $B$, acting on itself via left multiplication
and the projection $e_A$ coming from $\ve_A$, viewed as a linear operator
on $B$. It is well-known that $B_1$ is isomorphic as a vector space to
$Be_AB$ (here we identify $B$ with $\la(B)$, the algebra
of linear operators on $B$ coming from
left-multiplication by elements of $B$).
Moreover, the multiplication in $Be_AB$ is defined by
\begin{equation}\label{basicm}
(b_1e_Ab_2)(b_3e_Ab_4)=b_1\ve_A(b_2b_3)e_Ab_4.
\end{equation}

Assume now that the trivial representation $\one$ appears in
the second tensor power of the representation $V$ with multiplicity 1,
and it only appears in even tensor powers of $V$.
Decomposing $V^{\otimes n}$
as a direct sum of simple $\U=U_q\g$-modules, we define $V^{\otimes n}_{old}$
to be the direct sum of those simple modules which already appeared
in $V^{\otimes n-2}$. 
By semisimplicity and definition of  $V^{\otimes n}_{old}$,
we have a unique decomposition
$V^{\otimes n}=V^{\otimes n}_{old}\oplus V^{\otimes n}_{new}$
in these cases.
The following result has already more or less appeared before
in various publications; in the form below, 
see e.g. \cite{exc}, Prop. 4.10.

\begin{proposition}\label{oldstuff} Let $\Ca_n=\End_\U(V^{\otimes n})$.
Then the  algebra $\End_\U(V^{\otimes n+1}_{old})$ is isomorphic
to Jones' basic construction for $\Ca_{n-1}\subset \Ca_n$. In particular,
it is isomorphic as a vector space to $(\Ca_n\otimes 1)p_n(\Ca_n\otimes 1)$
where $p$ is the projection onto
$\one\subset V^{\otimes 2}$, and 
$p_n=1_{n-1}\otimes p\in \End_\U(V^{\otimes n+1})$.
\end{proposition}

We also need the following well-known property of  categorical
traces.

\begin{lemma}\label{markovproperty}
If $V, W$ are $\U$-modules with $V$ being irreducible, 
$a\in \End_\U(V^{\otimes 2})$ and $b\in \End_\U(W\otimes V)$,
then $tr((b\otimes 1_V)(1_W\otimes a))=tr(b)tr(a)$.
\end{lemma}

$Proof.$ This follows from the categorical definition of $tr$ and is
well-known. E.g. the proof of \cite{OW}, Prop. 1.4(c) can easily be
modified to prove the claim.

\section{Spinors for quantum groups}

\subsection{Roots and weights for orthogonal Lie groups} 
 For information
about roots and weights, see e.g. \cite{Hum},\cite{Kc}, and about
spinor representations, see e.g. \cite{Wy}. 
Let $\{ \e_i,\ 1\leq i\leq k\}$ be
the usual standard basis of $\R^k$. 
We represent the simple roots $(\alpha_i)$,
$i=1, 2, \ ...\ k$ of Lie types $B_k$ and $D_k$ as usual by
$\al_i=\e_i-\e_{i+1}$ for $i<k$, and as $\al_k=\e_k$ for Lie type $B_k$
and as $\al_k=\e_{k-1}+\e_k$ for Lie type $D_k$.
With these notations, the weight lattice is given by
$\Z^k \cup (\Z^k +\ve)$, where $\ve= (1/2,\ ...\ 1/2))$.
The irreducible representations of the corresponding Lie algebras
are labeled by the dominant weights $\la=(\la_i)$
which can be explicitly described as the set of all weights
$\la $ satisfying
$\la_1\geq\la_2\geq\ ...\geq \la_k\geq 0$ for type $B_k$
resp.
$\la_1\geq\la_2\geq\ ...\geq |\la_k|$ for type $D_k$.

\subsection{Pin groups}\label{pin}
The unique compact connected and simply connected Lie group
corresponding to the root systems $B_k$ and $D_k$ is the 
spin group $Spin(N)$ with $N=2k+1$ for type $B_k$ and
with $N=2k$ for type $D_k$. It is a 2-fold covering of
the orthogonal group $SO(N)$. As usual, we
denote the corresponding covering group of $O(N)$ by
$Pin(N)$. We embed $g\in O(2k)\to (g,\epsilon)\in O(2k)\times \Z/2
\subset SO(2k+1)$, where the sign is chosen so that
we obtain determinant one. This
embedding carries over to an embedding of $Pin(2k)$ into $Spin(2k+1)$,
i.e. we can consider  $Pin(2k)$ as a subgroup of  $Spin(2k+1)$.
As already indicated in the previous section, it will be convenient to
consider $Pin(2k)$ instead of $Spin(2k)$.

Algebraically, $Pin(2k)$ is a semidirect product of $Spin(2k)$ with $\Z/2$.
On the Lie algebra level, the  $\Z/2$-action is given by interchanging
the generators labeled by $k-1$ and $k$, i.e. the generators  belonging
to the endpoints of the $D_k$ graph next to its triple vertex.
This $\Z/2$-action induces a linear map
$\la\mapsto \bar\la$ on the weight space determined by permuting the
roots $\al_{k-1}$ and $\al_k$, and leaving the other simple roots fixed.
It is easy to check that if $\la=(\la_1,\ ...,\ \la_k)$,
then $\bar\la=(\la_1,\ ...,\la_{k-1},- \la_k)$.
The connection between irreducible $Spin(2k)$ and irreducible 
$Pin(2k)$-modules is described easily as follows:

- If $\la\neq \bar\la$ (i.e. $\la_k\neq 0$), then there exists a unique
irreducible $Pin(2k)$-module whose restriction to $Spin(2k)$ decomposes
as a direct sum of highest weight modules labeled by $\la$ and $\bar\la$.
We shall denote this  $Pin(2k)$-module
by $V_\la$ with $\la$ the dominant weight satisfying $\la_k>0$.

- If $\la =\bar\la$, there exist exactly two irreducible 
nonisomorphic $Pin(2k)$-modules, denoted  by $V_\la$ and $V_{\lad}$
whose restriction to $Spin(2k)$
is isomorphic to the highest weight module labeled by $\la$. 
Observe that in this case $\la$ can be
identified with a Young diagram and one takes for $\lad$ the Young
diagram with the same columns as $\la$ except that the first one now
has $2k-\la_1'$ boxes (where $\la_1'$ is the number of boxes in the first
column of $\la$). For all other dominant weights we define $\lad=\la$.

\subsection{Spinors}\label{spinors} 
 Let $S$ be the spinor module as constructed
via the Clifford algebra in Section \ref{classicalcase}.
In the odd-dimensional case, Lie type $B_k$ it is the irreducible representation
with highest weight $\La_k=\ve$, the fundamental weight dual to $\al_k$.
In the even-dimensional case, Lie type $D_k$, the module $S$ remains 
irreducible as a $Pin(2k)$-module, but decomposes into the direct sum
of two irreducible $Spin(2k)$-modules
whose highest weights are the fundamental weights $\La_{k-1}=\ve-\e_k$ and 
$\La_k=\ve$.

The module $S$ has the following properties:
Its  weights  are given by $\{ \om, \om= \frac{1}{2}
\sum_{i=1}^k \pm \e_i\}$, which holds for $S$ being viewed as a $Spin(2k+1)$
module as well as a $Pin(2k)$-module.
The following  tensor product rules for spinor groups
are well-known and follow easily from general
theory. More specialized treatments can also be found
in e.g. \cite{Ba}, \cite{BS} to name but a few.
If $V_\la$ is an irreducible module with highest weight $\la$
for Lie type $B_k$,
then 
$$
V_\la\otimes S\cong \bigoplus_\mu V_\mu,
$$
where the summation
goes over all dominant weights $\mu=\la+\om$ with $\om$ a weight
of $V$. For Lie type $D_k$, we have the following modification:
$$
V_\la\otimes S\ \cong\ V_\lad\otimes V\ \cong\ \bigoplus_\mu V_\mu 
\ \oplus\  \bigoplus_{\mu=\bar\mu} V_{\mud},
$$
where the summation
goes over all dominant weights $\mu=\la+\om$ with $\om$ a weight
of $S$ and with $\mu_k\geq 0$.

In the following we use Young diagram notation for labeling
the irreducible representations 
of the orthogonal groups, with the convention for $\lad$ as described above.
In particular we will write $[1^r]$ for the Young diagram with $r$
boxes in one column. We will need the following straightforward
examples, which are elevated to the rank of a lemma

\begin{lemma}\label{example} We have the following decompositions:
\begin{enumerate}
\item[(a)] $S^{\otimes 2}\cong \bigoplus_{s=0}^k V_{[1^{k-s}]}$ 
for Lie type $B_k$,
\item[(b)] $S^{\otimes 2}\cong \bigoplus_{s=-k}^{k} V_{[1^{k-s}]}$ 
for Lie type $D_k$,
\item[(c)] $S^{\otimes 3}\cong \bigoplus_{r=0}^k m_rV_{\e+[1^{k-s}]}$,
where the multiplicity $m_r$ is equal to $r+1$ for Lie type $B_k$,
and it is equal to $2r+1$ for Lie type $D_k$,
\item[(d)] For $Pin(2)\cong O(2)$ (see discussion in Section \ref{centralizers})
and $Spin(3)\cong SU(2)$, 
$S^{\otimes n}_{new}$ consists of one irreducible representation,
except for $n=2$ in the $Pin(2)$ case, where it is the direct sum
of two irreducible representations.
\end{enumerate}
\end{lemma}

\subsection{$q$-Dimensions}\label{qdimsect}  Recall that the $q$-dimension of a representation
is given by Eq. \ref{qdim}. We need more explicit formulas for certain
representations. We use the notation
$[n]_q= (q^{n}-q^{-n})/(q-q^{-1})$. 
Let $\U$ be equal to $U_qso_N$ (for $N$ odd) or the semidirect
product of  $U_qso_N$ with $\Z/2$ as in Section
\ref{centralizers} for $N$ even.
It is well-known that for
the $\U$ module $V_{[1^r]}$  with highest weight
$(1,\ ..., 1,0, ..., 0)$ (with $r$ 1s), we have 
\ignore{
$$\dim_qV_{[1^r]}\ =\ d(r,N)\ =\ \frac{[N-1]+[r]}{[r]} \prod_{j=2}^r 
\frac{[N-1-r+j]}{[r+1-j]}\ =:\ \binom{d}{r}_{qq}.$$}
\begin{equation}\label{dimforma}
\dim_qV_{[1^r]}\ =\ d(r,N)\ =\ \binom{N-1}{r}_{q}+\binom{N-1}{r-1}_{q}
= \binom{N}{r}_q\frac{q^{k-r}+q^{r-k}}{q^k+q^{-k}},
\end{equation}
where $\binom{n}{m}_{q}=[n]_q!/[m]_q![n-m]_q!$, and where
$[n]_q!=\prod_{i=1}^n [i]_q$.
This can be derived from the character formulas (see Eq. \ref{qdim})
or it can be read of as a special case from
the formulas in Section 5 of \cite{WBMW} for $\la=[1^r]$. 
Also observe that the dimension $\dim_qS=[2^k]_{qq}$ of the spinor module $S$ is given by
\begin{equation}\label{dimformc}
[2^k]_{qq} := 
\begin{cases} (q^{1/2}+q^{-1/2})(q^{3/2}+q^{-3/2})\ ...
\ (q^{k-1/2}+q^{1/2-k})  & \text{if $N=2k+1$ is odd,}\\
2(q^{}+q^{-1})(q^{2}+q^{-2})\ ...
\ (q^{k-1}+q^{1-k})
&\text{if $N=2k$ is even,}
\end{cases}
\end{equation}

\subsection{Symmetric representations of 
$\Ca_3$}\label{symmetric}
The calculation of the structure coefficients will be
significantly simplified by the existence of a certain involutive
antihomomorphism $^T$ on $\Ca_3$. In fact, it can be defined
for all $\Ca_n$, $n\in \N$. It satisfies
the involutive property $(c^T)^T=c$ for all $c\in\Ca_n$, $n\in\N$
and the functorial property
$$(c_1\otimes c_2)^T=c_1^T\otimes c_2^T\quad c_1\in \Ca_{n_1},\ c_2\in
\Ca_{n_2}.$$
The existence of this antihomomorphism is a consequence
of  Kashiwara's inner product on modules of
quantum groups (see e.g. \cite{exc}, Section 1.4 for details).

\begin{lemma}\label{lsymmetric}
For any simple module of $\End_\U(S^{\otimes 3}_{old})\subset\Ca_3$, 
we can find a basis $(v_r)$  of simultaneous
eigenvectors of $\Ca_2\otimes 1$ for which also each element in
$1\otimes \Ca_2$ is given by a symmetric matrix. In particular,
if $p^{(k)}$ is the projection onto $\1\subset S^{\otimes 2}$,
we can choose the basis such that $p^{(k)}_2=1\otimes p^{(k)}$ 
is given by the matrix 
$$(\sqrt{\dim_qV_{[1^i]}}\sqrt{\dim_qV_{[1^j]}}/\dim_qS)_{ij}.$$
\end{lemma}

\ignore{
In particular,
we can choose the basis such that the vector $\b$
is an eigenvector for all $a\in 1\otimes \Ca_2$,
where  $\b=(\sqrt{d(r,N)})v_r$,
with $d(r,N)$ as in \ref{dimforma}, and $0\leq r\leq N$ for $N$ even,
and $0\leq r\leq (N-1)/2$ for $N$ odd.}

$Proof.$
The first statement is a special case
of e.g. \cite{exc}, Lemma 1.9. In our special case of $\Ca_3$
path basis vectors would just be simultaneous
eigenvectors for the projections $p_1^{(s)}$, with $0\leq |s|\leq k$,
and the argument of that Lemma works for any element $A_2$ 
in $1\otimes \Ca_2$. 
The matrix of $p_2^{(k)}$ with respect to this basis
can be calculated by observing that it coincides with
the projection of Jones' basic construction for $\Ca_1\subset\Ca_2$
under the isomorphism stated in Prop. \ref{oldstuff}. 
Indeed, matrices for such projections have been calculated
in \cite{Su} or \cite{RaW} in terms of the weight vectors
of the trace; for the latter, see  Sections \ref{tracecont}
and \ref{qdimsect}. The second statement follows from this.

\ignore{
By Proposition \ref{oldstuff}, $p_2=(1\otimes p)$ in 
$\End_\U(S^{\otimes 3}_{old})$ acts as an orthogonal projection
onto $\C 1\cong \End_\U(S)\otimes 1\subset\End_\U(S^{\otimes 2})$
with respect to the orthonormal basis $(v_r)$, where $v_r=tr(e_r)^{-1/2}e_r$
and $e_r$ the projection onto $V_{[1^r]}\subset S^{\otimes 2}$.
Hence $1=\sum_r e_r=\sum_r tr(e_r)^{1/2}v_r$ is an eigenvector
of $1\otimes p$, and hence of any element $a$ of
the abelian algebra $1\otimes \Ca_2$. 
As $tr(e_r)=\dim_qV_{[1^r]}/\dim_q(S)^2$, the element 1 coincides
with the vector $\b$ in the statement up to a multiple.}

\subsection{Centralizers}\label{centralizers}
In the following we denote by $\U$ either
the Drinfeld-Jimbo $q$-deformation $U_q\so_{2k+1}$ of the
universal enveloping algebra of $\so_{2k+1}$
or the semidirect product of  $U_q\so_{2k}$ with $\Z/2$ as described in
Section \ref{pin}; the action described there
carries over to the quantum group level in the obvious way. 
So, in particular, the action of $t$ on a weight vector $v_\om$
of the spinor representation $S$ or the vector representation $V$
in the usual normalization is given by
\begin{equation}\label{taction}
t:\ v_\om\mapsto v_{\bar\om}.
\end{equation}

We define the algebras $\Ca_n^{(0)}=\C$ and  
$\Ca_n^{(N)}= \End_U(S^{\otimes n})$,
where $S$ is the spinor representation of type $B_k$ (for $N=2k+1$ odd)
or $D_N$ (for $N=2k$ even).
We will give an inductive procedure how to determine
the structure of these algebras. To get the induction started,
let us review the cases for $Pin(2)$ and $Spin(3)$, which are 
well-known:

$Pin(2)$: It is easy to check that $Pin(2)$ is isomorphic to the
orthogonal group $O(2)$: The two-fold covering of a circle is again
a circle, which defines an isomorphism between $Spin(2)$ and $SO(2)$.
This isomorphism extends to one between $Pin(2)$
and $O(2)$. Moreover, the spinor representation can be identified with
the usual two-dimensional representation of $O(2)$ under this isomorphism.
It is well-known that in the group case the centralizer algebra 
$\Ca_n^{(2)}$ is isomorphic to a quotient of Brauer's centralizer algebra
$D_n(2)$, whose simple components are labeled by Young diagrams whose
first two columns contain at most two boxes. It is also well-known that
we do not have a quantum deformation in this case, which was shown on
the categorical level in e.g. \cite{TbW2}, Lemma 7.5.
Hence we have in general that $\Ca_2^{(2)}\cong \C^3$ and
$\Ca_{n+1}^{(2)}$ is isomorphic to a direct sum of Jones' basic
construction for $\Ca_{n-1}^{(2)}\subset \Ca_{n}^{(2)}$ and
a one-dimensional direct summand labeled by the Young diagram $[n+1]$
(see \cite{Brauer}, \cite{Wbr}).

$Spin(3)\cong SU(2)$: Here the spinor representation $S$ corresponds to the
two-dimensional representation of $SU(2)$. The centralizer algebras are
well-understood in this case in the classical as well as in the quantum
case. They are given as quotients of Hecke algebras of type $A$, which
are also known as Temperley-Lieb algebras. Again, also in the
type $B$ case, $\Ca_{n+1}^{(3)}$ can be determined inductively as
a direct sum of Jones' basic
construction for $\Ca_{n-1}^{(3)}\subset \Ca_{n}^{(3)}$ and
a one-dimensional direct summand (see \cite{Jo}).

We will need the following notations for the general induction:
The element in $\Ca_2^{(N)}$ which projects onto the
submodule $V_{[1^{k-s}]}\subset V^{\otimes 2}$ will be denoted by $p^{(s)}$.
Observe that $p^{(k)}$ is the projection onto the trivial representation,
and that $s$ can be negative for type $D$ (see the end of
Section \ref{spinors}). Moreover, if $a\in \Ca_2^{(N)}$,
we define
$$a_i = 1_{i-1}\otimes a\otimes 1_{n-i-1}
\in\End_U(S^{\otimes n})= \Ca_n^{(N)}.$$

\begin{theorem}\label{main1}  The structure of $\Ca_n^{(N)}$,
as defined above, is determined for $n=2$ by Lemma \ref{example},
and it is determined inductively for $n>2$ and $N>0$ by
$$\Ca_{n+1}^{(N)}\ \cong\ \Ca_{n+1}^{(N-2)}\oplus \Ca_n^{(N)}p_n^{(N)}
\Ca_n^{(N)},$$
as a direct sum of algebras, with the multiplicative structure as in 
Eq \ref{basicm}.
In particular, $\Ca_n^{(N)}$  is generated by the elements 
$a_i$, $i=1,2,\ ...\ n-1$, with $a\in \End_\U(S^{\otimes 2})$.
\end{theorem}

$Proof$. We shall prove this theorem by induction on both $N$ and $n$.
Observe that the theorem follows from the discussion above for
$N=2$ and $N=3$.
Also observe that the trivial representation appears in the second
tensor power of the spinor representation $S$ in general for arbitrary $N$.
Hence, by Prop. \ref{oldstuff},
$\End_\U(S^{\otimes n+1}_{old})\cong \Ca_n^{(N)}p_n^{(N)} \Ca_n^{(N)}$.
So it suffices to show that $\End_\U(S^{\otimes n+1}_{new})\cong
\Ca_{n+1}^{(N-2)}$. The statement about the generators follows by 
induction on $n$ and $N$.

By Lemma \ref{prelim}(a), the commutant of the action of
the quantum group $\U'$ of type
$B_{k-1}$ resp. type $D_{k-1}$ on  $S^{\otimes n+1}[0]$ is isomorphic
to $\Ca_{n+1}^{(N-2)}$. Hence it suffices to show that the $\U$-module
generated by $V^{\otimes n+1}[0]$ is equal to $V^{\otimes n+1}_{new}$,
where now $\U$ is the  quantum group of type
$B_{k}$ resp. type $D_{k}$; indeed as
every irreducible $\U'$ submodule in $S^{\otimes n+1}[0]$ generates
an irreducible $\U$ module, we have 
$\Ca_{n+1}^{(N-2)}\cong\End_{\U'}(S^{\otimes n+1}[0])
\cong \End_\U(S^{\otimes n+1}_{new})$.

As $|\om_i|\leq 1/2$ for any weight $\om$ of $S$,
it follows from the tensor product rules for $S$
(see Section \ref{spinors}) by induction that also $|\la_1|\leq n/2$ for
any highest weight $\la$ in $S^{\otimes n}$. On the other hand,  it
is easy to check by induction, that any module $V_\la$ resp $V_\lad$
labeled by a dominant weight $\la$ with $|\la_1|\leq n/2$ does indeed
occur in $S^{\otimes n}$: the claim is obviously true for $n=1$,
and given a dominant weight $\la$ with $\la_1\geq 1$, we can always
find a weight $\om$ of $S$ such that $\la'=\la-\om$ is dominant.
As $V_{\la'}\subset V^{\otimes (n-1)}$ by induction assumption,
we obtain $V_\la\subset V_{\la'}\otimes S\subset S^{\otimes n}$.
Hence $S^{\otimes n}_{new}$ is a direct sum of highest weight modules
$V_\la$ such that $\la_1=n/2$, and any irreducible submodule 
of $S^{\otimes n}$ with such a highest weight is contained in  
$S^{\otimes n}_{new}$. 

As $n/2=n\langle \Lambda_1,\la\rangle$,
it follows that all these highest weight vectors are contained
in $S^{\otimes n}[0]$. Hence $S^{\otimes n}_{new}$ is contained in
the $\U$-module generated by  $S^{\otimes n}[0]$.
The other inclusion   follows from the fact that any highest weight
vector in  $V^{\otimes n}[0]$ has a weight $\la$ satisfying
$\lan \la,\La_1\ra =n/2$ (see the discussion before Lemma \ref{prelim}).
Hence it is contained in  $S^{\otimes n}_{new}$. This finishes the proof.

\begin{corollary}\label{cormain1} If $N=2k$ is even, the algebra
$\End_{Pin(2k)}(S^{\otimes l})$ is a quotient of $U(k,l)$, as defined in 
Def. \ref{algdef}.
If $N=2k+1$, the algebra $\End_{Spin(2k+1)}(S^{\otimes l})$ is a quotient
of $Uo(k,l)$.
\end{corollary}

$Proof.$ We have seen that the element $C$ has $N+1$ distinct eigenvalues,
and, for $N$ odd, the element $C^2$ has $(N+1)/2$ distinct eigenvalues.
Hence $C$ resp $C^2$ generates $\End_{Pin(N)}(S^{\otimes 2})$ for $N$ even,
resp $\End_{Spin(N)}(S^{\otimes 2})$ for $N$ odd. The rest follows from 
Theorem \ref{main1}.

\section{Structure Coefficients}

\subsection{Invariant elements, Clifford approach}\label{qCliffapp}
We have seen in Section \ref{classicalcase}
that the element $2C=\sum e_i\otimes e_i \in
Cl^{\otimes 2}\cong \End(S^{\otimes 2})$ generates
the commutant of the action of $Pin(2k)$, where $(e_i)_i$
is an orthonormal basis.
For the odd-dimensional case, it is convenient to consider the
restriction of the action on the module in the last sentence to
$Pin(N-1)$. We now denote this $Pin(N-1)$ module by $\tilde S$. 
It decomposes into a direct sum of irreducible
$Pin(N-1)$-modules $S_+\oplus S_-$, which are isomorphic as
$Spin(N-1)$-modules. Here we take as invariant element the canonical
element $C$ for the inner product of an $N-1$-dimensional subspace
of $V$. If we take the usual weight vectors as before, we can
express $C$ in the form
\begin{equation}\label{Csetup}
C(v_\la\otimes
v_\mu)=\frac{1}{2}\al_{\la,\mu}^{(N/2)}v_{\bar\la}\otimes
v_{\bar\mu} +\sum_{1\leq j<N/2,\ \la_j\neq\mu_j}
\al_{\la,\mu}^{(j)}v_{\bar\la^j}\otimes v_{\bar\mu^j};
\end{equation}
here $\bar\gamma^j$ is defined to coincide with $\gamma$ except for
a sign change in the $j$-th coordinate for any weight $\gamma$ of
$S$ resp. $\tilde S$ (here we use the notation for weights for
$Pin(N)$). In the even-dimensional case, the expression for $C$ is
as above, except that the sum goes until $N/2$ without the special
case for $j=N/2$. It is not hard to calculate the coefficients
$\al_{\la,\mu}^{(j)}$ which are equal to $\pm 1$.

In principle, this approach can be extended to the setting of
quantum groups, using a $q$-Clifford algebra, which has already
been studied (see \cite{Hay}, \cite{DF}). However, in this
context, one would also have to deform the multiplication of the
second tensor power of the $q$-Clifford algebra in a nontrivial way.
This makes calculations cumbersome. Instead, we shall produce the
$q$-analog of the invariant element $C$ by a
straightforward calculation of the coefficients
$\al_{\la,\mu}^{(j)}$ in the quantum case.

\subsection{Invariant elements, direct approach}
As motivated in the previous subsection, we now determine
a special element $C\in \End_\U(S^{\otimes 2})$ by
finding suitable coefficients for the expression
in Eq. \ref{Csetup}.

\ignore{
\begin{equation}\label{deffC}
C(v_\mu\otimes v_\nu)=\sum_{\mu,\nu} \sum_{j=1}^k
a_{\mu,\nu}^{j,-j} v_{\mu+\e_j}\otimes v_{\nu-\e_j} +
a_{\mu,\nu}^{-j,j}v_{\mu-\e_j}\otimes v_{\nu+\e_j}.
\end{equation}
}

\begin{proposition}\label{calc1}
Let $\U=U_qso_{2k}\rtimes \Z/2$. The element $C\in \End(S^{\otimes 2})$, defined by
$$C(v_\mu\otimes v_\nu)=\sum_j \delta_{\mu_j,-\nu_j}
(-q)^{\{\nu-\mu\}_{j-1}} v_{\bar\mu^j}\otimes  v_{\bar\nu^j}$$
commutes with the action of $\U$.
Here $\{\gamma\}_{j-1}=\sum_{i=1}^{j-1}\gamma_i$  for any $\gamma\in\R^k$.
Moreover $\bar\gamma^j$ is defined to coincide with $\gamma$ except for a
sign change in the $j$-th coordinate.
\end{proposition}

$Proof.$
It is easy to check that $C$ leaves invariant the weight spaces.
So it does commute with the generators $K_i$ of $\U$.
Also, as the action of $C$ on $v_\mu\otimes v_\nu$ only 
depends on the last coordinates of $\mu$ and $\nu$
as far as whether they are equal or not, $C$ also commutes with
the generator $t$ of $\Z/2$ in $\U$ (see Eq \ref{taction}).
It remains to check the equation
\begin{equation}\label{comparison}
C\Delta(X_j)(v_\mu\otimes v_\nu)=\Delta(X_j)C(v_\mu\otimes v_\nu),
\end{equation}
for $X_j=E_j, F_j$ and $1\leq j\leq k$. Recall that the coproduct is
defined by
$$\Delta(X_j)=K_j^{1/2}\otimes X_j+X_j\otimes K_j^{-1/2}.$$
It will be convenient to use the following notations
for the matrix coefficients of $C$:
\begin{equation}\label{coeff}
a_{\mu,\nu}^{(j)} = C_{\mu,\nu}^{\bar\mu^j,\bar\nu^j}
= \delta_{\mu_j,-\nu_j}(-q)^{\{\nu-\mu\}_{j-1}}.
\end{equation}
Comparing
the coefficients of the vector $v_{\mu +\e_j}\otimes v_{\nu-\e_{j+1}}$
in \ref{comparison} for $X_j=E_j$,  we obtain
\begin{equation}\label{comp3}
q^{-\langle \nu ,\al_j\rangle/2}a_{\mu+\al_j,\nu}^{(j+1)}+
q^{\langle \mu,\al_j\rangle/2}a_{\mu,\nu+\al_j}^{(j)}
=
q^{-\langle \nu -\e_{j+1},\al_j\rangle/2}a_{\mu,\nu}^{(j+1)}+
q^{\langle \mu +\e_j,\al_r\rangle/2}a_{\mu,\nu}^{(j)}.
\end{equation}
A similar equation follows if we consider the coefficients of the vector
$v_{\mu-\e_{j+1},\nu+\e_j}$.

Let us first consider the case $\U=U_qso_4$.  We can write the weights
and basis vectors
of $S$  as pairs of signs, e.g. $\ve=(++)$;
similarly, the basis vectors of $S\otimes S$ are written
as a vector with two such pairs, such as e.g. $(++,++)$. Then we deduce the
following identities from Eq \ref{comparison}:
\begin{align} \notag
a_{(+-,-+)}^{(j+1)}&=-q^{-1}a_{(-+,+-)}^{(j)},\quad   &a_{(-+,+-)}^{(j+1)}&=-q^{-1}a_{(+-,-+)}^{(j)},  \\
a_{(++,--)}^{(j+1)}&=-q^{-1}a_{(--,++)}^{(j)},\quad   &a_{(--,++)}^{(j+1)}&=-q^{-1}a_{(++,--)}^{(j)}, \label{equations} \\
a_{(+-,++)}^{(j+1)}&=a_{(-+,++)}^{(j)},\quad  & a_{(++,-+)}^{(j+1)}&=a_{(++,+-)}^{(j)}.
 \notag
\end{align}
Indeed, e.g. the equations in the first line follow from Eq. \ref{comparison} for the vector $(-+,-+)$
by comparing the coefficents of the vectors $(++,--)$ and $(--,++)$, see also Eq \ref{comp3}.
The other equations can be derived similarly, using generators $E_1,E_2, F_1, F_2$. Essentially,
these are calculations within $U_qsl_2$, applied to tensor products of vectors of weights $\pm 1$ or 0.

The general case with $\U=U_qso_{2k}$, with $k>2$
is not much more complicated. If we check
the claim in Eq. \ref{comparison} with $\langle \e_j,\al_r\rangle\neq 0$,
i.e. $j\in \{ j,j+1\}$ for $j<k$, only the $j$-th and $(j+1)$-st
coordinates of $\mu$ and $\nu$ are relevant for checking
Eq. \ref{comparison}, up to a common multiple for both sides.
We again get equations as in \ref{equations} from which we can determine the
coefficients $a_{\mu,\nu}^{(j)}$ by induction on $j$, starting with 
$a_{\mu,\nu}^{(1)}=\delta_{\mu_1,-\nu_1}$.

\begin{lemma}\label{eigenvector} (a) If $N=2k$ even, the eigenvalues of the
map $C$ are $[j]=(q^j-q^{-j})/(q-q^{-1})$ for  $-k\leq j\leq k$. 
The corresponding eigenspaces are $V_{[1^{2k-j}]}\subset S^{\otimes 2}$,
see Lemma \ref{example}.

(b) If $N=2k+1$ is odd, $C$ has the eigenvalues $[j+1/2]$ for $-k-1\leq j\leq k$.
\end{lemma}

$Proof.$ 
Let $\v=\sum_\la (-q)^{\langle \ve-\la,\rho\rangle}v_\la\otimes v_{-\la}$, where 
$\rho=(k-i)_i$ and $\ve$ is the highest weight vector of $S$.
 Then the $v_\la\otimes v_{-\la}$ coordinate of $C\v$ is given by 
$$\sum_j a_{\bar\la^j,-\bar\la^j}^{(j)} (-q)^{\langle \ve-\bar\la^j,\rho\rangle}
=(-q)^{\langle \ve-\la,\rho\rangle}\sum_j 
q^{\langle \bar\la^j-\la,\rho\rangle -\{ 2\la\}_{j-1}},$$
where we used the fact that $\langle \bar\la^j-\la,\rho\rangle +j-1$  is even.
Hence $\v$ is an eigenvector of $C$ if
 we can show that the set of exponents of $q$, namely
$\{ 2\la_j(k-j)-\{ 2\la\}_{j-1}, 1\leq j\leq k\}$ coincides with the set
of numbers $k+1-2r$, $1\leq r\leq k$. This is easily shown by induction
on $k$, by using the induction assumption for the weight
$\mu=(\la_2,\la_2,\ ....,\ \la_k)$ and observing that for $j=1$
we get $\pm (k-1)$ depending on the sign of $\la_1$.
This shows that $[k]$ is an eigenvalue.  Changing the sign of
the coefficient of $v_\la\otimes v_{-\la}$ for which $\la$ is a weight in $S_-$,
one also sees that $-[k]$ is an eigenvalue.

To prove the claim for the other values, observe that
$C$ leaves invariant the span $S^{\otimes 2}_r$ spanned by
vectors $v_\mu\otimes v_\mu$ for which $\mu_j=\nu_j=+$ if $j>r$.
Moreover, the action onto this subspace coincides with the one
of the element $C$ for $so_{2r}$. Hence we also have the eigenvalues
$\pm [r]$ for any $0\leq r<k$. As $C\in\End_\U(S^{\otimes 2})$ can have
at most $2k+1$ distinct eigenvalues, the claim follows. 

Part (b)
can be shown similarly.

\subsection{Odd-dimensional case} The same method also
works in the odd-dimensional case.  We shall do the case $O(3)$
in detail.  We shall consider a faithful representation of
$Cl(3)$ on a simple module of $Cl(4)$. We again use notation
$(++),\ (+-)$\ ... for the basis vectors. 
Then it is easy to see that the maps
$$E:\quad (--)\ \mapsto\  (++),\quad (-+)\ \mapsto\  (+-),$$
and $F$ being the transposed of $E$ with respect to this
basis define a representation of $U_qsl_2$.
It is the direct sum of two simple two-dimensional representations
with highest weight vectors $(++)$ and $(+-)$.
 Using the coproduct
as in the proof of Prop. \ref{calc1}, one can determine a
commuting operator $C$ as in Eq. \ref{Csetup}. If we set
 $\al_{\la,\mu}^{(1)}=1$ for any $\la,\mu$ with $\la_1=-\mu_1$, we
can determine
$\al^{(2)}_{\la,\mu}=1/[2]$ if both $v_\la$ and $v_\mu$
are highest weight vectors,
and
$$ \al^{(2)}_{(++,--)}=\al^{(2)}_{(+-,-+)}=
-q^{-1}/[2],\quad \al^{(2)}_{(--,++)}=\al^{(2)}_{(-+,+-)}=-q/[2],$$
where $[2]=q^{1/2}+q^{-1/2}$. Moreover, the coefficients for any tensor
product of two basis vectors are one of the above, where it only
depends whether the tensor factors are a highest or a lowest weight
vector.

\begin{proposition}\label{calc2}
If $N=2k+1$ is odd, we can determine coefficients $\al_{\la,\mu}^{(j)}$,
$1\leq j\leq k+1$
for $C$ as in Eq. \ref{Csetup} such that $C$ commutes with $\U=U_qso_N$
on $\tilde S^{\otimes 2}$, and that $C$ has the eigenvalues
$[j+1/2]=(q^{j+1/2}-q^{-j-1/2})/(q-q^{-1})$ for $-k-1\leq j \leq k$.
\end{proposition}

$Proof.$ If $j\leq k$, we take for $\al_{\la,\mu}^{(j)}$ the value as
in Prop. \ref{calc1}. For $\al^{(k+1)}_{\la,\mu}$, we define
 $\ve,\kappa$ to be  the 'vectors' consisting of the $k$-th and $(k+1)$-st
components of $\la$ and $\mu$ respectively.
If $\al^{(k)}_{\la,\mu}\neq 0$, we multiply it 
by $a^{(2)}_{\ve,\kappa}=-q^{\pm 1}/[2]$
as  in the $O(3)$-case to get $\al^{(k+1)}_{\la,\mu}$.
If $\al^{(k)}_{\la,\mu}= 0$, we set 
$\al^{(k+1)}_{\la,\mu}= \al^{(k)}_{\bar\la^k,\mu}/[2]$,
where $\bar\la^k$ coincides with $\la$ except for the $k$-th coordinate.
The claim about the coefficients of $C$ now
 follows from Prop. \ref{calc1} and the calculations for
the $O(3)$ case.

The claim about the eigenvalues is shown as in Lemma \ref{eigenvector},
where we now pick as eigenvector 
$\v=\sum (-q)^{\langle \ve-\lambda,\rho\rangle}v_\la\otimes v_{-\la}$, 
where $\rho=(k+1/2-i)_i$ and $\ve$ is the highest weight vector of $S$.

\subsection{Action in third tensor power} The main result of this
subsection is listed in Lemma \ref{shiftone}. It is elementary.
We will first deal with the slightly easier case $N$ even.
Recall that the $i$-th antisymmetrization $\bigwedge^i V$ of the
vector representation $V$ of $O(N)$ appears with multiplicity 1 in
$S^{\otimes 2}$ and that $\bigwedge^i V\otimes S$ contains a unique
summand which is isomorphic to $S$, which we denote by $S_i$, i.e.
we have a direct summand $S_i\subset S^{\otimes 3}$ defined by
$$S\cong S_i\subset \bigwedge^i V\otimes S\subset S^{\otimes 3},\quad 0\leq i\leq N.$$
Let now $(v_i)$ be an orthonormal basis of highest weight vectors
$v_i\in S_i$ of weight $\e$. Their span is a module of the commutant
of the $\U$ action. Let
\begin{equation}\label{veclincomb}
v_i=\sum \al^{(i)}_{\mu_1, \mu_2, \mu_3}
v_{\mu_1}\otimes  v_{\mu_2}\otimes v_{\mu_3},
\end{equation}
with $v_{\mu_j}$ a weight vector of $S$ for all indices $j$. We
extend the partial order of weights to tensor products of weight
vectors in alphabetic order, i.e. the order structure is
determined by the first factor for which the weights
are not the same.  Let $\Lambda_i$ be the highest weight of  $\bigwedge^i V$
for $i\neq N/2$ (in this case it is also irreducible as a $U_qso_N$
module). Then $v_\e\otimes v_{\Lambda_i-\e}\otimes v_{\e-\Lambda_i}$
and $v_{\bar\e}\otimes v_{\Lambda_i-\bar\e}\otimes v_{\e-\Lambda_i}$
are two maximal vectors with nonzero coefficients in the linear
combination of $v_i$ and $v_{N-i}$ for $i<N/2$. Hence $(1\otimes
C)v_i$ and $(1\otimes C)v_{N-i}$ are linear combinations of the
vectors $v_j$ with $j\leq i+1$ or $j\geq N-i-1$.

Moreover recall that $S=S_+\oplus S_-$ as a $U_qso_N$ module. Then
$\bigwedge^i V$ is contained in  $S_+^{\otimes 2}\oplus S_-^{\otimes
2}$ if $N/2-i$ is even, and in $S_+\otimes S_-\oplus S_-\otimes S_+$
if $N/2-i$ is odd. Hence we also have that $(1\otimes C)v_i$ is a
linear combination of vectors $v_j$ such that $j-i$ is odd. We have
set up everything for $N$ even for the following Lemma.

\begin{lemma}\label{shiftone} Let $C$ be the linear map
in $\End_\U(S^{\otimes 2})$ (for $N$ even) resp. in
$\End_\U(\tilde S^{\otimes 2})$ as defined in the previous sections.
The vector $(1\otimes C)v_i$ is a linear combination of $v_{i-1}$
and $v_{i+1}$.
\end{lemma}

$Proof.$ It follows from the definitions that $t$ has to map a highest weight
vector $v_\la$ of an $so_N$ module to a highest weight vector.
As  $\bigwedge^i V$ remains irreducible as an $so_N$ module for $i\neq N/2$
and as $t^2=1$, we have $tv_\la=\pm v_\la$ in this case.
One can now check directly for the highest weight vectors
of  $\bigwedge^i V$ that it is $+1$ for  $i<N/2$ and $-1$ for $i>N/2$. 
Hence we can normalize the
vectors $v_i$ and $v_{N-i}$ such that in their basis expansions the
vector $v_\e\otimes v_{\Lambda_i-\e}\otimes v_{\e-\Lambda_i}$ has
the same coefficient, and the coefficients of the vector
$v_{\bar\e}\otimes v_{\Lambda_i-\bar\e}\otimes v_{\e-\Lambda_i}$
differ by a sign. It follows from the definition of $C$ that
$(1\otimes C)v_i$ is a linear combination of the vectors $v_j$ with
$j\leq i+1$ and with $j-i$ odd, for $i<N/2$. A similar result also
holds for $(1\otimes C)v_{N-i}$.

Finally, we use Kashiwara's inner product (see e.g. \cite{Lu},  Section
\ref{symmetric} or also \cite{exc}, Section 1.4) for
representations of quantum groups. We actually only need the fact
that it is multiplicative for tensor products, and that its adjoint
maps $\End_{U_q\g}(W)$ into itself for any $U_q\g$ module $W$.
As the weight spaces are
mutually orthogonal and $\End_\U(S^{\otimes 2})$ is abelian, one
deduces that $C^T=C$, and hence also $(1\otimes C)$ is self-adjoint.
Hence, after suitably
normalizing the mutually orthogonal vectors $v_i$, we can assume $C$
to be symmetric. The claim for $N$ even follows from this and the statements at
the end of the last paragraph. The proof for $N$ odd goes exactly
the same way, after the notation is set up the right way.
This will be done in the remainder of this subsection.

\medskip

If $N$ is odd, we consider the module $\tilde S\cong \tilde S_+\oplus \tilde S_-$,
where $\SS_{\pm}\cong S$
with highest weight vectors $\ve_+=\ve$ and $\ve_-=\bar\ve$;
this is exactly the same decomposition of the corresponding $Pin(N+1)$-module
into a direct sum of irreducible $Spin(N+1)$-modules. 
Then we have as before that as $O(N)$-modules, $S_+^{\otimes 2}\cong
\oplus_{i=0}^k \bigwedge^{2i}V$ and
$S_+\otimes S_-\cong \oplus_{i=0}^k \bigwedge^{N-2i}V$.
If $(N+1)$ is divisible by 4, also $v_{-\ve}\in\SS_+$. In this case,
we define $v_{2i}$ to be a  higest weight vector of the unique module
$S\cong S_{2i}\subset  \bigwedge^{2i}V\otimes \SS_+\subset \SS_+^{\otimes 3}$,
and we define $v_{N-2i}$ to be a  higest weight vector of the unique module
$S\cong S_{N-2i}\subset  \bigwedge^{N-2i}V\otimes \SS_+\subset 
\SS_+\otimes \SS_-^{\otimes 2}$. If $(N+1)$ is not divisible by 4,
we define the vectors $v_{2i}\in\SS_+^{\otimes 2}\otimes \SS_-$ and
$v_{N-2i}\in \SS_+\otimes \SS_-\otimes \SS_+$. 
Using the fact that $C$ maps $\SS_+^{\otimes 2}$ to $\SS_-^{\otimes 2}$
and $\SS_+\otimes \SS_-$ to $\SS_-\otimes \SS_+$, it should now be no problem 
for the
reader to adapt the proof of Lemma \ref{shiftone} for the case $N$ odd.

\subsection{Technical lemma} Let $\{i\}=q^i+q^{-i}$ and define $b_i$
up to a sign by
$$b_i^2=\binom{N}{i}_q\frac{\{ N/2-i\} }{\{ N/2\} }.$$
Then we have the following lemma.

\begin{lemma} Let $A$ be a
symmetric  $(N+1)\times (N+1)$ matrix with $a_{ij}=0$ unless
$|i-j|=1$. Then the entries of $A$ are completely determined
by one eigenvalue $\la$ and its corresponding eigenvector $\b$
via Eq. \ref{aii}.
In particular, if $A$ has the eigenvalue $[N/2]$ with
eigenvector $\b = (b_i)$, where $b_i$ is as defined above, then
$$a_{i,i+1}^2=a_{i+1,i}^2=\frac{[i+1][N-i]}{\{ N/2-i\}\{ N/2-i-1\} }.$$
\end{lemma}

$Proof.$ It is straightforward to show by induction on $i$, using
the equation $A\b = \la\b$, that
\begin{equation}\label{aii}
a_{i,i+1}=\frac{\la}{b_ib_{i+1}}(b_i^2-b_{i-1}^2+b_{i-2}^2\ ...\ ).
\end{equation}
Now let $\la=[N/2]$ and let its eigenvector $\b$ be given as
above. Then we also have
$$b_i^2=\binom{N-1}{i}_q+\binom{N-1}{i-1}_q,$$
from which we get
$$a_{i,i+1}^2=\frac{[N/2]^2\binom{N-1}{i}_q^2}{b_i^2b_{i+1}^2}.$$
Using the definitions, it now is straightforward to show the claim.
\medskip

Let $H,E,K^{\pm 1}$ be the usual generators of the Drinfeld-Jimbo
quantum group $U_qsl_2$,
and let $V_N$ be its $(N+1)$-dimensional simple representation.
It can be defined via a basis $\{ e_o, e_1,\ ...\ e_N\}$
of eigenvectors of $K$ which satisfies
\begin{equation}\label{sl2rep}
E.e_r=([N-r+1][r])^{1/2}e_{r-1},
\quad K.e_r=q^{N-2r}e_r,\quad F.e_r=([N-r][r+1])^{1/2}e_{r+1}.
\end{equation}

\begin{corollary}\label{sl2connection}
The element $1\otimes C\in \End_\U(S^{\otimes 3})$
is given with respect to the basis $(v_r)$
in Lemma \ref{lsymmetric} by
$$A=(K^{1/2}+K^{-1/2})^{-1/2}(E+F)(K^{1/2}+K^{-1/2})^{-1/2},$$
where the elements on the right hand side stand for the matrices representing
these quantum group elements in $V_N$.
\end{corollary}

$Proof.$ Let $N$ be even.  It follows from Lemmas \ref{lsymmetric}, \ref{eigenvector} and
 \ref{shiftone} that the matrix
$A$ representing $1\otimes C$ with respect to the basis $(v_r)$
satisfies the conditions of the lemma;
for the statement about $\b$ observe that
$p_2^{(k)}$ is an eigenprojection of $C_2$ in the notation
of Lemma \ref{lsymmetric}.  For $N$ odd,
we get the appropriate eigenvectors for $C^2$,
restricted to the basis vectors in $\SS_+^{\otimes 3}$, as well as
for the eigenvectors in $\SS_+\otimes \SS_-^{\otimes 2}$, by 
Lemma \ref{lsymmetric}. The claim can now be shown also for $C$,
using its block structure with respect to our basis.

\subsection{Relations for centralizer algebras}  We will need the following
nonstandard $q$-defor-

\noindent
mation of the universal enveloping algebra of $so_N$.
It was  defined by Gavrilik and Klimyk (see e.g. \cite{GK}) and by Noumi and
Sugitani \cite{NS}.  It is also 
a special case of the co-ideal subalgebras of $U_qsl_N$ defined by Letzter
for $\Theta = id$, see\cite{Le}, Remark 2.4. It is not isomorphic to the usual
Drinfeld-Jimbo quantum group, see  \cite{Le}, Remark 2.3.

\begin{definition}\label{coidealdef} (a) The algebra $U_q'so_l$ is defined via
generators $B_1,B_2,\ ...\ B_{l-1}$ and relations $B_iB_j=B_jB_i$ for $|i-j|>1$
and
$$B_i^2B_{i\pm 1}-(q+q^{-1})B_iB_{i\pm 1}B_i +B_{i\pm 1}B_i^2= B_{i\pm 1}.$$

(b) The algebra $U_q(l,k)$ is the quotient of $U'_qso_{l}$ defined via the 
additional relation $\prod_{j=-k}^k (B_i-[j])=0$.

(c) The algebra $Uo_q(l,k)$ is the quotient of the subalgebra of  $U_q'so_l$
generated by $B_i^2$, $1\leq i<l$ defined via the additional relation
$\prod_{j=1}^k (B_i^2-[j-1/2]^2)=0$.
\end{definition}

It is clear that for $q=1$ we obtain the relations for the universal enveloping algebra
of the Lie algebra $so_l$, see e.g. the remarks before Lemma \ref{solrelations}.

\begin{theorem}\label{main2} (a) If $N=2k$ is even, we have representations of
$\U=U_qso_N\rtimes \Z/2$ and $U_q'so_l$ on $S^{\otimes l}$ which are each 
others commutant. Moreover, the image of $U_q'so_l$ factors through $U_q(l,k)$.

(b) If $N=2k+1$ is odd, we have actions of $\U=U_qso_N$ and $Uo_q(l,k)$ on
$S^{\otimes l}$ which are each others commutant.
\end{theorem}

$Proof.$ In case (a), it follows from Lemma \ref{eigenvector} that $C$ has $N+1$
eigenvalues, and hence generates $\End_\U(S^{\otimes 2})$. It follows from
Theorem \ref{main1} that $\End_\U(S^{\otimes l})$ is generated by the elements $C_i$,
$1\leq i<l$.

It can easily be seen that the claim will follow if the commutation relations between
$C_1$ and $C_2$ are checked, using   Lemma \ref{eigenvector}.
Now observe that by Lemma \ref{eigenvector} and by Cor.  \ref{sl2connection}
the elements $C_1$ and $C_2$ act on the span of the vectors $(v_r)\subset S^{\otimes 3}$
via the matrices representing the elements  $D_1=(k^{1/2}-k^{-1/2})/(q-q^{-1})$
and $A$ in $U_qsl_2$. To check the commutation relations between $C_1$ and $C_2$,
observe that
$$[k,E]=(q^2-1)Ek=(1-q^{-2})kE,\quad [k^{-1},E]=(1-q^{2})k^{-1}E=(q^{-2}-1)Ek^{-1}.$$
The relations involving $F$ are obtained from above by substituting $E$ by $F$ and $q$ by $q^{-1}$.
Setting $D_1=(k^{1/2}-k^{-1/2})/(q-q^{-1})$, we obtain
$$[D_1, E]=q^{1/2}Ek^{1/2}+q^{-1/2}Ek^{-1/2}=E(q^{1/2}k^{1/2}+q^{-1/2}k^{-1/2}),$$
and a similar expression with $E$ replaced by $F$. Hence
\begin{align} \notag
[D_1,[D_1, E+F]]\ =\ & E(q^{1/2}k^{1/2}+q^{-1/2}k^{-1/2})^2+F(q^{-1/2}k^{1/2}+q^{1/2}k^{-1/2})^2\\
\ =\ &(q^{1/2}+q^{-1/2})^2(E+F)+(k^{1/2}-k^{-1/2})(E+F)(k^{1/2}-k^{-1/2}). \notag
\end{align}
This can be rewritten as
$$D_1^2(E+F)-(q+q^{-1})D_1(E+F)D_1+(E+F)D_1^2=(E+F).$$
As $(k^{1/2}+k^{-1/2})^{-1/2}$ commutes with $D_1$, and obviously also
$C_i$ commutes with $C_j$ provided $|i-j|>1$, we have proved the relation
in Def. \ref{coidealdef} for $i$ and $i+1$ for $i=1$. Proving the relation with
$B_1$ and $B_2$ interchanged can again be done via a calculation in $U_qsl_2$,
or one checks that one gets the same matrices in our set-up
if we  interchange $C_1$ and $C_2$. This shows that the relation holds
for the summand $\Ca_2^{(N)}p_2^{(N)}\Ca_2^{(N)}$ of $\Ca_3^{(N)}$
in Theorem \ref{main1}. The general claim follows from that theorem and induction
by $N$, observing that the projection $p^{(N)}_n$ is an eigenprojection of $C_n$.

Again, the proof for the case $N$ odd goes along the same lines.

\begin{remark}\label{generalization} We have stated our results here only for
$q$ generic over the ground field of the complex numbers. It is not hard to
generalize them to more general fields; e.g. even though square roots
appear in our proofs, we do not expect them to be essential. Moreover the
definition of the algebra $U_q'so_l$ only involves elements in the ring
$\Z[q,q^{-1}]$. So a similar result should also hold over that ring;
indeed, by general results of Lusztig's also centralizer algebras 
of suitably defined quantum groups
have canonical bases defined over $\Z[q,q^{-1}]$, see \cite{Lu}, Section 27.3.
\end{remark}

\begin{remark}\label{comparison} In previous work \cite{KW}, \cite{TbW}
tensor categories were classified whose Grothendieck semiring was the one
of a unitary, orthogonal or symplectic group. There an intrinsic description 
of centralizer algebras via braid group representations played a crucial
role. For spinor representations, it is more difficult to describe
these braid representations as the standard generators have too many
eigenvalues. This made it necessary to consider a new description which 
generalizes the braid relations. Indeed, in many ways, the algebra
$U_q'so_l$ can be considered an algebraic object of type $A_{l-1}$.
They should be useful in proving similar classification results for
spinor groups.
\end{remark}

\begin{remark}\label{Rowell} In the recent publication \cite{RW}
Rowell and Wang study certain representations of braid groups which
they call Gaussian representations. They are of particular interest
in their studies motivated by quantum computing. They conjecture that these
representations are
related to the centralizer algebras of quantum groups for
spinor representations for certain roots of unity. Again, as already
noted in Remark \ref{comparison}, it is diffcult
to determine these braid representations because of the increasing
number of eigenvalues. The detailed analysis of representations of $\Ca_3$
in this paper could be useful in solving this problem, see e.g. the proof of
Theorem \ref{main2}.
\end{remark}

\begin{remark}\label{Lehrer} It is not hard to deduce
from Theorem \ref{main1} that for $N$ odd the $R$-matrices
generate $\End_\U(S^{\otimes l})$ for all $l$. Indeed, using
Drinfeld's quantum Casimir, one shows that the $R$ matrix
has $(N+1)/2$ distinct eigenvalues, and hence generates
$\End(S^{\otimes 2})$. 
Similar results have been shown by
 Lehrer and Zhang for 
vector representations of classical type and certain representations
in type $A_1$ and $G_2$, see \cite{LZ}.
\end{remark}

\end{document}